\documentclass[a4paper,12pt]{amsart}

\usepackage{amsmath,amsthm,amssymb,amsfonts,amscd,graphicx,mathtools,color}
\usepackage[all]{xy} 
\numberwithin{equation}{section}
\theoremstyle{plain} 
	\newtheorem{thm}{Theorem}[section]
	
	\newtheorem{lem}[thm]{Lemma}
	\newtheorem{conj}[thm]{Conjecture}

	\newtheorem*{thm*}{Theorem}
	\newtheorem*{conj*}{Conjecture}
\theoremstyle{definition}
        \newtheorem{prop}[thm]{Proposition}
	\newtheorem{defn}[thm]{Definition}

	\newtheorem{exmp}[thm]{Example}
	
\theoremstyle{remark}
	\newtheorem{rem}[thm]{Remark}
        
	\newtheorem*{pf}{Proof}
\def\EE{{\mathbb E}}
\def\CC{{\mathbb C}}

\def\HH{{\mathbb H}}

\def\RR{{\mathbb R}}

\def\ZZ{{\mathbb Z}}

\def\D{{\mathcal D}}

\def\F{{\mathcal F}}

\def\I{{\mathcal I}}

\def\O{{\mathcal O}}
\def\P{{\mathcal P}}

\def\T{{\mathcal T}}
\def\U{{\mathcal U}}

\def\Z{{\mathcal Z}}
\def\p{{\partial }}
\def\-Mod{\text{-}{\rm Mod}}
\def\Mod{{\rm Mod}}
\def\Mod-{{\rm Mod}\text{-}}
\def\Pre-Tr{{\rm Pre}\text{-}{\rm Tr}}

\def\Hom{{\rm Hom}}

\def\det{{\rm{det}}}

\def\Stab{{\rm{Stab}}}

\def\End{{\rm End}}

\def\ds{/\hspace{-1.0mm}/}

\def\h{{\mathfrak h}}

\def\p{\partial }
\newcommand{\Res}{\mathop{\rm Res}}

\def\ns{{\nabla}\hspace{-1.4mm}\raisebox{0.3mm}{\text{\footnotesize{\bf /}}}}

\begin{document}
\title{A Frobenius manifold for $\ell$-Kronecker quiver}

\author{Akishi Ikeda} 
\address{Department of Mathematics, Josai University, Saitama, Japan}
\email{akishi@josai.ac.jp}
\author{Takumi Otani} 
\address{Department of Mathematics, Graduate School of Science, Osaka University, 
Toyonaka Osaka, 560-0043, Japan}
\email{u930458f@ecs.osaka-u.ac.jp}
\author{Yuuki Shiraishi}
\address{School of Science IUPS, Osaka University, 
Toyonaka Osaka, 560-0043, Japan}
\email{yshiraishi@iups.sci.osaka-u.ac.jp}
\author{Atsushi Takahashi}
\address{Department of Mathematics, Graduate School of Science, Osaka University, 
Toyonaka Osaka, 560-0043, Japan}
\email{takahashi@math.sci.osaka-u.ac.jp}

\maketitle
\thispagestyle{empty}
\begin{abstract}
We construct a Frobenius structure whose intersection form coincides with the generalized Cartan matrix of 
the $\ell$-Kronecker quiver $K_{\ell}$ and underlying complex manifold is isomorphic to the space of stability conditions for 
the bounded derived category of finitely generated modules over the path algebra $\CC K_{\ell}$. 
\end{abstract}
\section{Introduction}

A Frobenius manifold is a complex manifold whose tangent bundle is a flat family of Frobenius algebras 
with grading operator called the Euler vector field satisfying certain conditions (see Definition~\ref{def of Frobenius}). 
The notion of a Frobenius manifold is formulated by B. Dubrovin \cite{dub} in his study of 
differential equations and integrable systems (e.g., Painlev\'{e} equations) related to $2$-dimensional topological field theories.  
However, this structure (except for the potentiality) was firstly discovered by K. Saito and his collaborators \cite{sai:2,sys} 
in their study of the invariant ring of a finite reflection group and its orbit space. Its regular subspace coincides with 
the domain for period mappings of a certain differential form on cycles in the Milnor fiber of the corresponding singularity.
Later, inspired by the above construction, a Frobenius structure on the deformation space of 
an isolated hyper-surface singularity is constructed via the choice of a special differential form called a {\it primitive form} $\zeta$
in the de Rham cohomology group of the total space twisted by the differential of the singularity (see \cite{sai}).
The existence of a primitive form is proved by M. Saito \cite{msai} (for tame holomorphic functions by Douai--Sabbah \cite{ds:1, ds:2}).
In particular, the Gelfand--Leray form $\check{\zeta}$ of a primitive form gives the differential form for 
the period mappings mentioned above.

T.~Bridgeland introduced the notion of a space of stability conditions for a triangulated category 
(constructed from a symplectic manifold, a tame function with isolated critical points etc).
Roughly speaking, this space is a complex moduli manifold consisting of tuples $(Z,\P)$ of a group homomorphism
called a {\it central charge} $Z$ from the Grothendieck group to $\CC$ and an $\RR$-graded family of full additive categories $\P$
satisfying certain conditions (\cite{bri}).
The importance of the space of stability conditions is that this complex manifold is conjectured to be the moduli space
of deformations of a mathematical object mentioned inside the bracket above (e.g. see \cite{tak,kst:1,kst:2}). 
Indeed, this conjecture is verified for some cases concerning the derived categories of (resp. Calabi--Yau completions of) 
Fukaya--Seidel categories for tame functions with isolated critical points, equivalently, 
the derived categories of modules over the mirror path algebras. 
Namely, in \cite{bqs, hkk} (resp. \cite{ike,ike17,wang}), their central charges are given by the oscillatory integrals 
(resp. period integrals of Gelfand--Leray forms) for primitive forms and the spaces of stability conditions are isomorphic to 
the (resp. universal cover of regular) orbit spaces of corresponding Weyl group invariant theories. 
More precisely, the following commutative diagrams summarize their results for $A_{2}$-singularity $F_{{\bf t}}(z)=z^{3}+t^{2}z+t^{1}$: 
\begin{table}[htbp]\label{table:1}
  \begin{center}
    \begin{tabular}{cc}
\begin{minipage}{0.5\hsize}
\begin{center}
\[
\xymatrix{
\h/W \ar[rr]^{\cong}  \ar[d]_{\Pi} & & \Stab(\D^{b}(A_{2})) \ar[d]^{\Z} \\
T_{{\bf p}}(\h/W)\cong \h \ar[rr]^{\cong} & & \Hom(K_{0}(\D^{b}(A_{2}), \CC),
}
\]
\end{center}
\end{minipage}
\quad
\begin{minipage}{0.5\hsize}
\begin{center}
\[
\begin{dcases}
\displaystyle
\Pi({\bf t})=\int e^{F_{\mathbf{t}}(z)}\zeta\in T_{{\bf p}}(\h/W),\\ 
\h^{*}_{\ZZ}\cong H_{1}(\CC, {\rm Re}(F_{{\bf p}}(z))\ll 0; \ZZ),\\
\zeta=dz.
\end{dcases}
\]
\end{center}
\end{minipage}
\end{tabular}
  \end{center}
    \end{table}

\begin{table}[htbp]\label{table:2}
  \begin{center}
    \begin{tabular}{cc}
\begin{minipage}{0.5\hsize}
\begin{center}
\[
\xymatrix{
\widetilde{\h_{\rm reg}/W} \ar[rr]^{\cong}  \ar[d]_{\Pi_{N}} & & \Stab(\D_{N}(A_{2})) \ar[d]^{\Z_{N}} \\
T_{{\bf p}}(\widetilde{\h_{\rm reg}/W})\cong \h \ar[rr]^{\cong} & & \Hom(K_{0}(\D_{N}(A_{2})), \CC),
}
\]
\end{center}
\end{minipage}
\quad
\begin{minipage}{0.5\hsize}
\begin{center}
\[
\begin{dcases}
\displaystyle
\Pi_{N}({\bf t})=\int \check{\zeta}\in T_{{\bf p}}(\widetilde{\h^{\rm reg}/W}),\\ 
\h^{*}_{\ZZ}\cong H_{N}((\Sigma_{N}F_{{\bf p}})^{-1}(0); \ZZ) , \\
\displaystyle
\check{\zeta}=\Res \frac{dx_{1}\wedge \cdots dx_{N}\wedge dz}{\Sigma_{N}F_{{\bf t}}}.\\
\end{dcases}
\]
\end{center}
\end{minipage}
\end{tabular}
  \end{center}
    \end{table}
\noindent       
Here we denote by $\h$ (resp. $\h_{\rm reg}$) the Cartan subalgebra (resp. without reflection hyper-planes), 
by $\D^{b}(A_{2})$ (resp. $\D_{N}(A_{2})$) the bounded derived category of finitely generated modules over the path algebra $\CC A_{2}$ 
(resp. of finite total dimension dg-modules over the Calabi--Yau $N$-completion of $\CC A_{2}$),
by $\Z$ (resp. $\Z_{N}$) the local homeomorphism $(Z, \P) \mapsto Z$, by $\Sigma_{N}F_{{\bf t}}$ the $N$-th suspension 
$\displaystyle F_{{\bf t}}+\sum_{i=1}^{N}x^{2}_{i}$ and the tilde stands for the universal cover. 
Note that the integral dual lattice $\h^{*}_{\ZZ}$ is naturally isomorphic to the relative homology group of Lefschetz thimbles 
(resp. the middle dimensional homology group of vanishing cycles) and integrations are taken after analytic continuations 
from the reference point ${\bf p}$ to ${\bf t}$ (see also \cite[Section 5 (5.7)]{sai}).    

Frobenius structures on the orbit spaces of the Weyl groups are constructed for finite Coxeter groups \cite{sys, dub:99},
extended affine Weyl groups \cite{dz, dszz}, elliptic Weyl groups and Jacobi groups \cite{ber:1, ber:2, sat}
(see also \cite{kamase,kms} for Saito structures without metrics). 
The purpose of the present paper is to construct the Frobenius structure from the invariant theory of 
the Weyl group associated to $\ell$-Kronecker quiver $K_{\ell}$:
\begin{thm}[Theorem \ref{thm : Frobenius manifold for l--Kronecker quiver}]\label{intro:main}
Let $W$ be the Weyl group of $K_\ell$ and $\widetilde{X}$ be the universal cover of a certain subspace $X$ of $\h$ (see Definition \ref{defn : definition of X}).
There exists a unique Frobenius structure of rank $2$ and dimension $1-\dfrac{2}{h}$ on $\widetilde{X}\ds W$ 
such that $\displaystyle e=\frac{\p}{\p t^{1}}, \ E=t^{1}\dfrac{\p}{\p t^{1}}+\dfrac{2}{h}t^{2}\dfrac{\p}{\p t^{2}}$
for the flat coordinates $(t^{1}, t^{2})$ in Subsection \ref{flat coords 4.1} and 
the intersection form coincides with the generalized Cartan matrix $A$ of $\ell$-Kronecker quiver $K_{\ell}$. 
\end{thm}
Here $\displaystyle h:=\dfrac{2\pi\sqrt{-1}}{\log\rho}$ where $\rho>1$ is the eigenvalue of the Coxeter transformation. 
Since $\displaystyle \rho=\exp\left(\dfrac{2\pi\sqrt{-1}}{h}\right)$ tautologically, 
the number $h$ can be regard as a generalization of the Coxeter number. 
Inspired by the conjectural relation between Frobenius structure and spaces of  stability conditions explained above,
the covering space of regular orbit subspace $\widetilde{X}_{\rm reg} / W (\subset \widetilde{X}\ds W)$ is obtained as  
the space of stability conditions for the bounded derived category of nilpotent modules over the preprojective algebra 
$\Pi_2 (K_{\ell})$ associated to $K_{\ell}$. 
Moreover, due to \cite{dk}, the space $\widetilde{X}\ds W$ is isomorphic to the space of stability conditions for 
the bounded derived category of finitely generated modules over the path algebra $\CC K_{\ell}$ 
(see Proposition \ref{prop:space of X ds W}). 
We expect stronger results that the deformed flat coordinates on the Frobnius manifold $\widetilde{X} \ds W$ coincide with the central charges on $\Stab(\D^{b}(K_{\ell}))$ 
(see Conjecture \ref{conj : Frobenius manifold and satbility conditions for generalized Kronecker quiver}) .


\bigskip
\noindent
{\it Acknowledgement.}
A. I is supported by JSPS KAKENHI Grant Number 16K17588 and partially by JSPS KAKENHI Grant Number 16H06337.
Y. S is supported by JSPS KAKENHI Grant Number 19K14531.
A. T is supported by JSPS KAKENHI Grant Number 16H06337.

\section{Frobenius manifolds and periods}
We recall the definition of the Frobenius manifold and related properties.  
\subsection{Frobenius manifolds}

The original definition and their basic properties are given by B.~Dubrovin \cite{dub}. 
For the notations, we use the following definition:
\begin{defn}[{\cite[Definition 2.1]{st}}]\label{def of Frobenius}
Let $M=(M,\O_{M})$ be a connected complex manifold or a formal manifold over $\CC$ of dimension $\mu$
whose holomorphic tangent sheaf and cotangent sheaf are denoted by $\T_{M}, \Omega_M^1$ respectively
and $d$ a complex number.
A {\it Frobenius structure of
rank $\mu$ and dimension $d$ on M} is a tuple $(\eta, \circ , e,E)$, where $\eta$ is a non-degenerate $\O_{M}$-symmetric bilinear
form on $\T_{M}$, $\circ $ is $\O_{M}$-bilinear product on $\T_{M}$, defining an associative and commutative
$\O_{M}$-algebra structure with the unit $e$, and $E$ is a holomorphic vector field on $M$, called
the Euler vector field, which are subject to the following axioms:
\begin{enumerate}
\item The product $\circ$ is self-adjoint with respect to $\eta$: that is,
\begin{equation*}
\eta(\delta\circ\delta',\delta'')=\eta(\delta,\delta'\circ\delta''),\quad
\delta,\delta',\delta''\in\T_M. 
\end{equation*} 
\item The {\rm Levi}--{\rm Civita} connection $\ns:\T_M\otimes_{\O_M}\T_M\to\T_M$ with respect to $\eta$ is
flat: that is, 
\begin{equation*}
[\ns_\delta,\ns_{\delta'}]=\ns_{[\delta,\delta']},\quad \delta,\delta'\in\T_M.
\end{equation*}
\item The tensor $C:\T_M\otimes_{\O_M}\T_M\to \T_M$  defined by 
$C_\delta\delta':=\delta\circ\delta'$, $(\delta,\delta'\in\T_M)$ is flat: that is,
 \begin{equation*}
\ns C=0.
\end{equation*} 
\item The unit element $e$ of the $\circ $-algebra is a $\ns$-flat homolophic vector field: that is,
\begin{equation*}
\ns e=0.
\end{equation*} 
\item The bilinear form $\eta$ and the product $\circ$ are homogeneous of degree $2-d$ ($d\in\CC$), $1$ respectively 
with respect to Lie derivative ${\rm Lie}_{E}$ of {\rm Euler} vector field $E$: that is,
\begin{equation*}
{\rm Lie}_{E}(\eta)=(2-d)\eta,\quad {\rm Lie}_{E}(\circ)=\circ.
\end{equation*}
\end{enumerate}
\end{defn}
Next we expose some basic properties of the Frobenius manifold without their proofs.
Let us consider the space of horizontal sections of the connection $\ns$:
\[
\T_M^f:=\{\delta\in\T_M~|~\ns_{\delta'}\delta=0\text{ for all }\delta'\in\T_M\}
\]
which is a local system of rank $\mu $ on $M$ such that the metric $\eta$ 
takes constant value on $\T_M^f$. Namely, we have 
\begin{equation*}
\eta (\delta,\delta')\in\CC,\quad  \delta,\delta' \in \T_M^f.
\end{equation*}
\begin{prop}[{\cite[Definition 2.2]{st}}]\label{prop:flat coordinates}
At each point of $M$, there exist local coordinates $(t^{1},\dots,t^{\mu})$, called flat coordinates, such that
$e=\p_1$, $\T_M^f$ is spanned by $\p_1,\dots, \p_{\mu}$ and $\eta(\p_i,\p_j)\in\CC$ for all $i,j=1,\dots, \mu$,
where we denote $\p/\p t^{i}$ by $\p_i$. 
\end{prop} 
The axiom $\ns C=0$, implies the following:
\begin{prop}[{\cite[Proposition 2.4]{st}}]\label{prop:potential}
At each point of $M$, there exist the local holomorphic function $\F$, called Frobenius potential, satisfying
\begin{equation*}
\eta(\p_i\circ\p_j,\p_k)=\eta(\p_i,\p_j\circ\p_k)=\p_i\p_j\p_k \F,
\quad i,j,k=1,\dots,\mu,
\end{equation*}
for any system of flat coordinates. In particular, one has
\begin{equation*}
\eta_{ij}:=\eta(\p_i,\p_j)=\p_1\p_i\p_j \F. 
\end{equation*}
\end{prop}
Denote by $\eta^{ij}$ the $(i,j)$-entry of the matrix $(\eta_{ij})^{-1}$.

\begin{exmp}[{\cite[Example 1.1]{dub}}]\label{important example}
Let $M$ be a Frobenius  manifold of rank $2$ and dimension $d$ whose flat coordinates are $(t^{1},t^{2})$. 
If $d\ne-1,1,3$, then the Frobenius potential $\F$ of $M$ is given as follows:
\begin{equation*}
\label{Frobeniuspotentialofgeneraltype}
\F(t^1,t^2)= \frac{1}{2}\eta_{12}(t^1)^2t^2+c(t^2)^{\frac{3-d}{1-d}},
\end{equation*}
where $\eta_{12}\in\CC\backslash\{ 0 \}, ~c\in\CC$. 
\end{exmp}

\subsection{The intersection form and the second structure connection}
Let $M$ be a Frobenius manifold of rank $\mu$ and dimension $d$. 
We recall the intersection form of the  Frobenius manifold $M$.
\begin{defn}[{\cite[Appendix G (G.1)]{dub}}]
Define $\Delta$ as the determinant of the $\O_{M}$-endomorphism $C_{E}\in\End_{\O_{M}}(\T_{M})$: 
\begin{equation*}
\Delta:=\det(C_E)\in\O_{M},
\end{equation*}
called the {\it discriminant}. 
Set $\D:=\{ t\in M~|~\Delta(t)=0\}$ (called the discriminant locus) and $M_{\rm reg}:=M\backslash\D$
(called the regular subspace). 
\end{defn}

The intersection form of the Frobenius manifold $M$ is defined as follows: 
\begin{defn}[{\cite[Lecture 3 (3.13)]{dub}}]
\label{definitionofintersectionform}
Let $(s^1,\dots,s^\mu)$ be local coordinates of $M$. 
Set the $\O_{M}$--symmetric bilinear form $g:\Omega_{M}^{1}\times\Omega_{M}^{1}\to\O_{M}$ as 
\begin{equation*}
g(\omega,\omega'):=\sum_{a,b=1}^\mu\left\langle\frac{\p}{\p s^a},\omega\right\rangle\cdot\eta(ds^a,ds^b)
\cdot\left\langle C_E\frac{\p}{\p s^b},\omega'\right\rangle ,
\end{equation*}
where $\langle -,- \rangle$ is the contraction. 
The $\O_{M}$-bilinear form $g$ is called the {\it intersection form of the Frobenius manifold} $M$ .
\end{defn}
The intersection form $g$ is independent of a choice of local coordinates. 
The intersection form with respect to the flat coordinates $(t^{1},\dots,t^{\mu})$ is given by 
\begin{equation*}
g=\sum_{a,b=1}^{\mu}\eta^{ia}\eta^{jb}E\dfrac{\p^{2} \F}{\p t^{a}\p t^{b}}.
\end{equation*}
Denote by $\nabla$ the Levi--Civita connection of $g$. 
The connection $\nabla$ is often called the {\it Second structure connection} or the {\it Dubrovin connection}. 

\begin{prop}[{\cite[Theorem 9.4]{her}}]
The second structure connection $\nabla$ is a flat connection. 
\end{prop}

\begin{defn}[{\cite[Definition G.1]{dub}}]
Actions of the fundamental group $\pi_{1}(M_{\rm reg})$ on the universal covering of $M_{\rm reg}$
can be lifted to isometries in ${\rm Aut}(\CC^{\mu}, I)$ on $\CC^{\mu}\cong T^{*}_{p} M_{\rm reg}$ 
with respect to the bilinear form $I$ induced by the intersection form $g$. 
Denote by $\Theta$ the associated representation: 
\begin{equation*}
\Theta: \pi_{1}(M_{\rm reg}, p)\longrightarrow {\rm Aut}(\CC^{\mu}, I).
\end{equation*}
The group
\begin{equation*}
W_{M}:={\rm Im}\,\Theta
\end{equation*}
is called the {\it monodromy group of the Frobenius manifold}. 
\end{defn}
\begin{defn}[{\cite[Definition 6]{dub:2}}]
A function $x\in\O_{M_{\rm reg}}$ is called a {\it period} if the $1$-form $dx\in\Omega_{M_{\rm reg}}^{1}$ is flat with respect to 
the second structure connection $\nabla$, i.e., 
\begin{equation*}
\nabla dx =0 .
\end{equation*}
\end{defn}

\begin{prop}\label{degree two coordinate}
\label{prop:flat function}
Let $(x^{1},\dots,x^{\mu})$ be local flat coordinates on $M_{\rm reg}$ with respect to $g$. 
If $d\ne 1$, then the following function: 
\begin{equation*}
t(x^{1},\dots,x^{\mu}):=\sum_{i=1}^{\mu}g\left( \dfrac{\p}{\p x^{i}},\dfrac{\p}{\p x^{j}} \right) x^{i} x^{j}
\end{equation*}
is a flat function, namely, $\ns d t=0$. 
\end{prop}
\begin{pf}
The statement follows from similar arguments in \cite[Section 3.3, 3)]{sai} and \cite[Section 5.3)]{sai}.
\qed
\end{pf}

The following Lemma \ref{lem : dubrovin} will be the key lemma which enables us to reconstruct 
structure coefficients of the product $\circ$ from contravariant components of the Levi--Civita connection for 
the intersection form $g$
(see Subsection \ref{reconst of product from g}): 
\begin{lem}[{\cite[Lemma 3.4]{dub}}]
\label{lem : dubrovin}
Let $g$ be the intersection form of a Frobenius structure $M$ of rank $\mu$ and dimension $d$. 
The following equality holds: 
\begin{align*}
g(dt^{i},\nabla_{\frac{\p}{\p t^{k}}}dt^{j})=\left(\dfrac{d-1}{2}+d_{j}\right)\sum_{a,b=1}^{\mu}\eta^{ia}\eta^{jb}\p_k\p_a\p_b\F,
\end{align*} 
where $d_{j}$ is the degree of the flat coordinate $t_{j}$ with respect to $E$, i.e., 
$\displaystyle \ns_{\frac{\p}{\p t^{j}}} E=d_{j}\cdot \frac{\p}{\p t^{j}}$.
\end{lem}

\section{The generalized root system associated to $\ell$-Kronecker quiver}
Throughout this paper, we assume that $\ell\in \ZZ$ and $\ell\ge 3$.
We summarize results for the generalized root system associated to $\ell$-Kronecker quiver.
\begin{defn}
The $\ell$-Kronecker quiver $K_{\ell}$ consists of vertices $\{ 1,2\}$ and $\ell$ edges from $1$ to $2$:
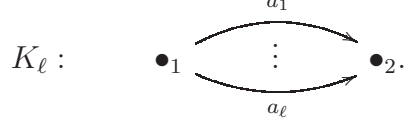
\[
\xymatrix{
K_{\ell} : & \bullet_{1} \ar@/_1pc/[rr]_{a_\ell} \ar@/^1.2pc/[rr]^{a_1} & \vdots & \bullet_{2}. 
}
\]
\begin{enumerate}
\item The {\it generalized Cartan matrix} associated to $\ell$-Kronecker quiver is defined as
\begin{equation*}
A:=
\begin{pmatrix}
a_{11} & a_{12} \\
a_{21} & a_{22}
\end{pmatrix}
=
\begin{pmatrix}
2 & -\ell \\
-\ell & 2
\end{pmatrix},
\end{equation*}
where $a_{ij}:=2\delta_{ij}-(q_{ij}+q_{ji})$, $\delta_{ij}$ is the Kronecker delta and $q_{ij}$ is the number of arrows 
from $i$ to $j$.
\item The {\it root lattice} $L$ associated to $K_{\ell}$ is the following free abelian group
\begin{equation*}
L:=\ZZ\alpha_{1}\oplus\ZZ\alpha_{2}
\end{equation*}
with generators $\alpha_{1}$ and $\alpha_{2}$, called {\it simple roots}. 
\item Define the symmetric bilinear form as follows:
\begin{equation*}
I:L\times L\to\ZZ; \quad I(\alpha_{i},\alpha_{j}):=a_{ij}.
\end{equation*}
Define the {\it reflection} with respect to the simple root $\alpha_{i}$ as the following isometric homomorphism $r_{i}:L\to L$
with respect to $I$:  
\begin{equation*}
r_{i} (\lambda):=\lambda-I(\lambda,\alpha_{i})\alpha_{i} , \quad \lambda\in L.
\end{equation*} 
\item The group $W:=\langle r_{1},r_{2} \rangle \subset {\rm Aut}(L, I)$ is called the {\it Weyl group}. 
Define the set of {\it real roots} $\Delta^{\rm re}$ as
\begin{equation*}
\Delta^{\rm re}:= \{w(\alpha_i)\in L~|~w\in W ,~i=1,2\}, 
\end{equation*}
\item Define the {\it Coxeter transformation} $c$ as $c:=r_{1}r_{2}\in W$. 
\end{enumerate}
\end{defn}

Let $\D^{b}(K_{\ell})$ be the bounded derived category of finitely generated modules over the path algebra associated to $K_{\ell}$, 
$\chi$ the Euler form on the Grothendieck group $K_{0}(\D^{b}(K_{\ell}))$, 
$S_{i}$ the simple module corresponding to the vertex $i$, 
$\Delta^{re}(\D^{b}(K_{\ell})):=\langle r_{[S_{1}]},r_{[S_{2}]} \rangle\{[S_{1}],[S_{2}]\}$
where $r_{[S_{i}]}(\lambda):=\lambda-(\chi+\chi^{T})(\lambda, [S_{i}])[S_{i}]$ for $\lambda, [S_{i}]\in K_{0}(\D^{b}(K_{\ell}))$
and $[S[1]]$ the automorphism on $K_{0}(\D^{b}(K_{\ell}))$ induced by the Serre functor $S$ shfited by $[1]$.  	
\begin{lem}
The tuples $(L, I, \Delta^{\rm re}, c)$ and $(K_{0}(\D^{b}(K_{\ell})), \chi+\chi^{T}, \Delta^{re}(\D^{b}(K_{\ell})), [S[1]])$ form
generalized root systems in the sense of \cite[Definition 2.1]{stw}.
Moreover, these tuples coincide with each other under the identification $[S_{i}]=\alpha_{i} \ (i=1,2)$.
\end{lem}
\begin{pf}
For the former assertion, see \cite[Section 2.2]{stw}. The latter assertion follows from straightforward calculations.
\qed
\end{pf}

The matrix representations of $r_{i}$ and $c$ (denote them by same symbols with $r_{i}$ and $c$) 
with respect to the basis $(\alpha_{1},\alpha_{2})$ are given by 
\begin{equation*}
r_{1}=
\begin{pmatrix}
-1&\ell\\
0&1
\end{pmatrix}
,\quad r_{2}=
\begin{pmatrix}
1&0\\
\ell&-1
\end{pmatrix}
,\quad c=
\begin{pmatrix}
\ell^2-1&-\ell\\
\ell&-1
\end{pmatrix}.
\end{equation*}
Let $\rho,\rho^{-1}$ be eigenvalues (spectral radius) of the Coxeter transformation $c$:
\begin{equation*}
\rho:=\dfrac{\ell^{2}-2+\sqrt{\ell^{4}-4\ell^{2}}}{2},\quad \rho^{-1}=\dfrac{\ell^{2}-2-\sqrt{\ell^{4}-4\ell^{2}}}{2} .
\end{equation*}

\subsection{The space $X$}
\begin{defn}
Let the group homomorphism $\alpha_{i}^{*}:L\to\CC$ be the dual of $\alpha_{i}$, i.e., 
$\langle\alpha_{i}^{*},\alpha_{j}\rangle:= \alpha_{i}^{*}(\alpha_{j})=\delta_{ij}$. 
Set
\begin{align*}
\h_{\RR}&:=\Hom_{\ZZ}(L,\RR)\cong\RR\alpha_{1}^{*}\oplus\RR\alpha_{2}^{*} \\
\h&:=\Hom_{\ZZ}(L,\CC)\cong\CC\alpha_{1}^{*}\oplus\CC\alpha_{2}^{*},
\end{align*}
and call $\h$ the {\it Cartan subalgebra}. 
Set the dual spaces 
\begin{align*}
\h_{\RR}^{*}&:= L\otimes_{\ZZ}\RR\cong\RR\alpha_{1}\oplus\RR\alpha_{2} \\
\h^{*}&:= L\otimes_{\ZZ}\CC\cong\CC\alpha_{1}\oplus\CC\alpha_{2} 
\end{align*}
The $W$-action on $L$ can be extended linearly to $\h^{*}$ and $\h^{*}_{\RR}$. 
Define the $W$-action on $\h$ and $\h_{\RR}$ by
\begin{equation*}
\langle w(Z),\lambda\rangle:=\langle  Z,w^{-1}(\lambda)\rangle\quad\text{for $Z\in\h$ and $\lambda\in\h^{*}_{\CC}$} .
\end{equation*}
\end{defn}
The vector spaces $\h,\h_{\RR},\h^{*},\h_{\RR}^{*}$ are equipped with the Euclidean topology for finite dimensional ones. 

\begin{defn}[cf. {\cite[section 5.1]{kac}}]
\label{roots}
Set 
\begin{equation*}
L_{+}:=\sum^{2}_{i=1}\ZZ_{\ge0}\alpha_{i}, \quad L_{-}:=-L_{+}=\sum^{2}_{i=1}\ZZ_{\le0}\alpha_{i}. 
\end{equation*}
\begin{enumerate}
\item Define the set of {\it positive real roots} $\Delta^{\rm re}_{+}$ and the one of {\it negative real roots} $\Delta^{\rm re}_{-}$ as 
\begin{equation*}
\Delta^{\rm re}_{+}:=\Delta^{\rm re}\cap L_{+},\quad \Delta^{\rm re}_{-}:=\Delta^{\rm re}\cap L_{-}.
\end{equation*}
\item Let $K :=\{ \lambda\in L_+\backslash\{0\} ~|~ I(\lambda,\alpha_i)\le 0,~ i=1,2\}$.
The set of {\it imaginary roots} $\Delta^{\rm im}$ is defined as $\Delta^{\rm im}:=\Delta^{\rm im}_{+}\cup\Delta^{\rm im}_{-}$,
where the set of {\it positive imaginary roots} $\Delta^{\rm im}_{+}$ and the one of {\it negative imaginary roots} $\Delta^{\rm im}_{-}$ 
are defined as 
\begin{equation*}
\Delta^{\rm im}_{+} :=W(K)=\{w(\lambda)~|~w\in W,~ \lambda\in K\},\quad\Delta^{\rm im}_{-}:=-\Delta^{\rm im}_{+}.
\end{equation*}
\end{enumerate}
\end{defn}

\begin{defn}
The {\it imaginary cone} $\I$ is defined as the closure of the convex hull of $\Delta^{\rm im}_{+}\cup\{0\}$:
\begin{equation*}
\I :=\overline{{\rm Conv}(\Delta^{\rm im}_{+}\cup\{0\})}\subset\h_{\RR}^{*}.
\end{equation*}
Set $\I_{0}:=\I\backslash\{0\}$ and call it the {\it blunt imaginary cone}.
\end{defn}
Let $\displaystyle \nu:=\dfrac{\ell+\sqrt{\ell^{2}-4}}{2}$ and $\displaystyle \nu^{-1}=\dfrac{\ell-\sqrt{\ell^{2}-4}}{2}$.
Obviously $\nu^{2}=\rho$.
We have
\begin{align*}
\I &=\left\{ z_{1}\alpha_{1}+z_{2}\alpha_{2}\in\h_{\RR}^{*}~\middle|~z_{1},z_{2}\ge0,~\sum_{i,j=1}^{2}z_{i}z_{j}I(\alpha_{i},\alpha_{j})\le0\right\} \\
&=\{ z_{1}\alpha_{1}+z_{2}\alpha_{2}\in\h_{\RR}^{*}~|~z_{1},z_{2}\ge0,~2(z_{1}^{2}-\ell z_{1}z_{2}+z_{2}^{2})\le0\} \\
&=\{ z_{1}\alpha_{1}+z_{2}\alpha_{2}\in\h_{\RR}^{*}~|~z_{1},z_{2}\ge0,~(z_{1}-\nu z_{2})(z_{1}-\nu^{-1}z_{2})\le0\} .
\end{align*}
\begin{defn}[{\cite[Definition 2.7]{ike}}]
\label{defn : definition of X}
For $\lambda \in \h^*$,  set $H_{\lambda}:=\{ Z\in\h~|~Z(\lambda)=0 \}$. 
Define $X\subset\h$ as
\begin{equation*}
X:=\h\backslash\bigcup_{\lambda\in \I_{0}}H_{\lambda} ,
\end{equation*}
and define $X_{\rm reg}\subset X$ as
\begin{equation*}
X_{\rm reg}:=X\backslash\bigcup_{\alpha\in\Delta^{\rm re}_+}H_\alpha .
\end{equation*}
\end{defn}

\begin{lem}[{\cite[Lemma 2.9]{ike}}] 
The sets $X$ and $X_{\rm reg}$ are open subsets of $\h$.
\end{lem}

\begin{prop}[{\cite[Proposition 2.17]{ike}}]\label{proper disconti and free}
The $W$-action can be restricted to $X$ and $X_{\rm reg}$. 
This $W$-action is properly discontinuous on $X$, in particular, is free on $X_{\rm reg}$.
\end{prop}

Let $\D^{b}(\Pi_2 (K_{\ell}))$ be the bounded derived category of nilpotent modules over the preprojective algebra 
$\Pi_2 (K_{\ell})$ associated to $K_{\ell}$, 
$\ZZ[2] \subset {\rm Aut}\left(\D^{b}(\Pi_2 (K_{\ell}))\right)$ the subgroup generated by the shift functor 
$[2]$.
Denote by ${\rm Br}(\D^{b}(\Pi_2 (K_{\ell})))$ the subgroup of ${\rm Aut}(\D^{b}(\Pi_2 (K_{\ell})))$ generated by 
spherical twists. 
\begin{thm}[{\cite[Theorem 1.1]{ike}}]
There is a covering map
\[
\underline{\pi}: {\rm Stab}^{\circ}\left(\D^{b}(\Pi_2 (K_{\ell}))\right) \longrightarrow X_{\mathrm{reg}} / W
\]
and the subgroup 
$\ZZ[2] \times {\rm Br}\left(\D^{b}(\Pi_2 (K_{\ell}))\right) \subset {\rm Aut}\left(\D^{b} (\Pi_2 (K_{\ell}))\right)$ 
acts as the group of deck transformations.
\end{thm}

\begin{lem}\label{coordinate change}
Define the new basis $(\beta_{1},\beta_{2}):=(\alpha_{1},\alpha_{2})P$ where
\begin{equation*}
\displaystyle
P:= 
\begin{pmatrix}
\frac{\nu}{\sqrt{(\ell^{2} -4)\nu}} & \frac{1}{\sqrt{(\ell^{2} -4)\nu}} \\
\frac{1}{\sqrt{(\ell^{2} -4)\nu}} & \frac{\nu}{\sqrt{(\ell^{2} -4)\nu}}
\end{pmatrix}.
\end{equation*}
\begin{enumerate}
\item The matrix representations $R_{i}$ of reflections $r_{i}$ and the one of $I$ 
on $\h^{*}_{\CC}$ with respect to the basis $(\beta_{1},\beta_{2})$ are given by
\begin{equation*}
R_{1}:=P^{-1} r_1 P =
\begin{pmatrix}
0 & \nu \\
\nu^{-1} & 0
\end{pmatrix}
,\quad
R_{2}:=P^{-1} r_2 P =
\begin{pmatrix}
0 & \nu^{-1} \\
\nu & 0
\end{pmatrix}
,\quad
P^T A P =
\begin{pmatrix}
0 & -1 \\
-1 & 0
\end{pmatrix}.
\end{equation*}
We also have $P^{-1}cP=P^{-1}r_{1}r_{2}P=R_{1}R_{2}$ and use the same symbol $c$ for $P^{-1}cP$.
\item Let $(\beta_{1}^{*},\beta_{2}^{*})$ be the dual basis for $(\beta_{1},\beta_{2})$.
Then $(\beta_{1}^{*},\beta_{2}^{*})=(\alpha_{1}^{*},\alpha_{2}^{*})(P^{T})^{-1}$.
\item The matrix representations of reflections $r_{i} \ (i=1,2)$ on $\h$ 
with respect to the basis $(\beta_{1}^{*},\beta_{2}^{*})$
are given by $(R_{1}^{T})^{-1}$ and $(R_{2}^{T})^{-1}$ respectively.
\end{enumerate}
Denote by $(x^1,x^2)$ and $(x_1,x_2)$ the linear coordinates with respect to the basis 
$(\beta_{1}^{*},\beta_{2}^{*})$ of the Cartan subalgebra $\h$ 
and the basis $(\beta_1,\beta_2)$ of the dual $\CC$-vector space $\h^{*}_{\CC}$ respectively. 
\begin{enumerate}\setcounter{enumi}{3}
\item We have
\begin{equation*}
R_{1}\cdot (x^{1},x^{2})=(\nu^{-1} x^{2},\nu x^{1}),\quad R_{2}\cdot (x^{1},x^{2})=(\nu x^{2},\nu^{-1} x^{1}) \ .
\end{equation*}
\item We have
\begin{equation*}
\I=\{(x_1,x_2)\in\RR^2~|~x_1,x_2\ge0,~x_1x_2\le0\}, 
\end{equation*}
and hence 
\begin{align*}
X&=\CC^2\backslash\bigcup_{0\le a\le1}\{ (x^1,x^2)\in\CC^2~|~ax^1+(1-a)x^2=0\} \\
&\label{firstX} =\CC^2\backslash\bigcup_{0\le\lambda\le\infty}\{ (x^1,x^2)\in\CC^2~|~x^1=-\lambda x^2\}.
\end{align*}
\end{enumerate}
\end{lem}
\begin{pf}
The statements {\rm (i)} to {\rm (iv)} follow from straightforward calculations. 
The statement {\rm (v)} follows from {\cite[Lemma 2.8]{ike}}.
\qed
\end{pf}

\begin{lem}\label{fundamental group of X}
Define the cycle $\gamma:[0,1]\to X$ as 
\[
\gamma(t):=\left(e^{2\pi\sqrt{-1}t},e^{2\pi\sqrt{-1}t}\right)\in X.
\]
Then $\gamma$ is a generator of the fundamental group $\pi_{1}(X)$ of $X$ and hence $\pi_{1}(X)=\langle \gamma \rangle \cong\ZZ$. 
\end{lem}
\begin{pf}
By Lemma~\ref{coordinate change} {\rm (v)}, there exists an isomorphism of complex manifolds:
\begin{equation*}
\displaystyle X \xrightarrow{\cong} (\CC\setminus \{0\})\times (\CC\setminus \RR_{\le 0});
\quad (x^{1}, x^{2})\mapsto \left(x^{1}, \frac{x^{1}}{x^{2}}\right).
\end{equation*}
The cycle $\{(e^{2\pi\sqrt{-1}t},1)~|~t\in[0,1]\}\subset (\CC\setminus \{0\})\times (\CC\setminus \RR_{\le 0})$ is
an generator of the fundamental group 
$\pi_{1}((\CC\setminus \{0\})\times (\CC\setminus \RR_{\le 0}))\cong 
\pi_{1}(\CC\setminus \{0\})\times \pi_{1}(\CC\setminus \RR_{\le 0})\cong \ZZ$
and the image of $\gamma$ by the above isomorphism.
\qed
\end{pf}

\subsection{Weyl group invariant theory}
\label{subsection : Weyl group invariant theory}
Set 
\begin{equation}
\label{defn:universal covering  of X}
\widetilde{X}:=\left\{ \left( y^{1},y^{2} \right) \in\CC^{2}~\middle|~\left| {\rm Im}\left(y^{1}-y^{2}\right) \right|<\pi \right\}\subset\CC^{2} .
\end{equation}

\begin{prop}
\label{prop:universal covering}
The space $\widetilde{X}$ is the universal covering space of $X$. The covering map is given by
\begin{equation*}
\pi_{X}:\widetilde{X}\longrightarrow X,\quad (y^{1},y^{2})\mapsto (e^{y^{1}},e^{y^{2}})
\end{equation*}

\end{prop}
\begin{pf}
For $(y^{1},y^{2})\in\widetilde{X}$, set $\lambda:=e^{y^{1}-y^{2}}=e^{{\rm Re}\, (y^{1}-y^{2})}e^{\sqrt{-1}\,{\rm Im}\, (y^{1}-y^{2})}$. 
Since $\lambda$ is not a negative real number since $\left| {\rm Im}(y^{1}-y^{2}) \right|<\pi$ 
and $e^{y^{1}}=e^{y^{1}-y^{2}+y^{2}}=\lambda e^{y^{2}}$, the map $\pi_{X}$ is well-defined.
The statement follows from Lemma \ref{fundamental group of X}.
\qed
\end{pf}

\begin{rem}
Fix the isomorphism $\h\cong \CC^{2}$ via the basis $(\beta_{1}^{*},\beta_{2}^{*})$.
Define the map $\pi_{\CC^{2}}:\CC^{2}\rightarrow \CC^{2}\cong \h$ as $(y^{1},y^{2})\mapsto (e^{y^{1}},e^{y^{2}})$.
The boundary $\p \widetilde{X}$ of $\widetilde{X}$ is given by 
\[
\p \widetilde{X}=\left\{ (y^{1},y^{2})\in\CC^{2}~\middle|~\left| {\rm Im}\left(y^{1}-y^{2}\right)  \right|=\pi \right\}.
\]
We have $\pi_{\CC^{2}}(\p \widetilde{X})=\h\backslash X=\bigcup_{\lambda\in\I_{0}}H_{\lambda}$.
\end{rem}

\begin{rem}
If we choose the following domain as $\widetilde{X}$:
\begin{equation*} 
\left\{ \left( y^{1},y^{2} \right) \in\CC^{2}~\middle|~(2n-1)\pi<{\rm Im}\left(y^{1}-y^{2}\right)<(2n+1)\pi \right\} \quad (n\ne 0),
\end{equation*}
the monodromy group of the resulting Frobenius manifold does not coincide with the Weyl group for $K_{\ell}$. 
See also Proposition~\ref{monodromy prop}.
\end{rem}

We can lift the Weyl group action on $\widetilde{X}$ as follows: 
\begin{align}
\label{action:R1}
R_{1}\cdot (y^{1},y^{2}) &=(y^{2}-\log \nu,~y^{1}+\log\nu ) \\
\label{action:R2}
R_{2}\cdot (y^{1},y^{2}) &=(y^{2}+\log \nu,~y^{1}-\log\nu ).
\end{align}
The $W$-action above is equivariant to $\pi_{X}$.
The fundamental group $\pi_{1}(X)$ acts on the unversal covering sapce $\widetilde{X}$. 
This $\pi_{1}(X)$-action is given by
\begin{equation*}
\gamma\cdot (y^{1},y^{2})=(y^{1}+2\pi\sqrt{-1},~y^{2}+2\pi\sqrt{-1})
\end{equation*}
for the generator $\gamma\in\pi_{1}(X)$. 
These two actions on $\widetilde{X}$ can be extended naturally on $\CC^{2}(\supset\widetilde{X})$. 
Note that these actions on $\CC^{2}$ are properly discontinuous.

\begin{lem}
Set
\begin{equation*}
\widetilde{\U}:=\left\{ \left( y^{1},y^{2}\right)\in\widetilde{X} ~\middle|~ \left| {\rm Re}\left(y^{1}-y^{2}\right) \right|\le \log \nu \right\}\subset\widetilde{X},
\end{equation*}
\[
\U:=\pi_{X}\left(\widetilde{\U}\right)= \left\{ \left( x^{1},x^{2} \right) \in X ~\middle|~ \nu^{-1} \le \left| \dfrac{x^{1}}{x^{2}} \right| \le \nu \right\}\subset X.
\]
Then the subset $\widetilde{\U}$ is a fundamental domain of the $W$-action in \eqref{action:R1} and \eqref{action:R2} on $\widetilde{X}$.
The subset $\U$ is a fundamental domain of the $W$-action in Lemma \ref{coordinate change} {\rm (iv)} on $X$. 
\end{lem}

\begin{pf}
The statement follows from explicit presentations of actions \eqref{action:R1} and \eqref{action:R2}. 
\qed
\end{pf}

\begin{defn}\label{analytic spaces}
Let $W\curvearrowright \CC^{2}\supset \widetilde{X}$ be the $W$-action in \eqref{action:R1} and \eqref{action:R2}. 
Define the complex analytic space $\CC^{2}\ds W$ as follow: 
\begin{itemize}
\item Its underlying space is $\CC^{2}/W$, the quotient space of $\CC^{2}$ by $W$.
\item Let $\pi:\CC^{2}\to \CC^{2}/W$ be the quotient map.
Denote by $\O_{\CC^{2}}^{W}$ the $W$--invariant subsheaf of $\O_{\CC^{2}}$
the sheaf of germs of holomorphic functions for $\CC^{2}$. 
Define the structure sheaf $\O_{\CC^{2}\ds W}$ as $\O_{\CC^{2}\ds W}:=\pi_{*}\O_{\CC^{2}}^{W}$. 
\end{itemize}
Define the complex analytic space $\widetilde{X}\ds W$ as follows:
\begin{itemize}
\item Its underlying space is $\widetilde{X}/W$, the quotient space of $\widetilde{X}$ by $W$.
\item Let $\pi:\widetilde{X}\to \widetilde{X}/W$ be the quotient map.
Denote by $\O_{\widetilde{X}}^{W}$ the $W$--invariant subsheaf of $\O_{\widetilde{X}}$,
the sheaf of germs of holomorphic functions for $\widetilde{X}$. 
Define the structure sheaf $\O_{\widetilde{X}\ds W}$ as $\O_{\widetilde{X}\ds W}:=\pi_{*}\O_{\widetilde{X}}^{W}$. 
\end{itemize}
\end{defn}

\begin{prop}
\label{prop:invariant theory of CC}
Define the map $\varphi:\CC^{2}\to\CC^{2}; (y^{1},y^{2})\mapsto (s^{1},s^{2})$ where
\begin{align*}
s^{1}&= e^{\frac{h}{2}(y^{1}-y^{2})}-e^{\frac{h}{2}(y^{2}-y^{1})},\\
s^{2}&= y^{1}+y^{2}.
\end{align*}
Then the map $\varphi$ induces an isomorphism between complex analytic spaces 
\begin{equation}
\label{eqn : isomorphism from C2 ds W to C2}
\underline{\varphi}:\CC^{2}\ds W \overset{\cong}{\longrightarrow} \CC^{2},\quad [(y^{1},y^{2})]\mapsto (s^{1},s^{2}) .
\end{equation}
In particular, the complex analytic space $\CC^{2}\ds W$ is a complex manifold.
\end{prop}

\begin{pf}
We have $\underline{\varphi}\circ\pi=\varphi$ and  
$\underline{\varphi}_{*}\O_{\CC^{2}\ds W}=\underline{\varphi}_{*}\pi_{*}\O_{\CC^{2}}^{W}\cong\varphi_{*}\O_{\CC^{2}}^{W}$. 
We shall show that $\varphi_{p}:\O_{\CC^{2},\varphi(p)}\to\O_{\CC^{2},p}^{W}$ is isomorphic for any $p\in \CC^{2}$. 
This morphism $\varphi_{p}$ is given by 
\begin{equation*}
\displaystyle \varphi_{p}(f)(y^{1},y^{2})=f\circ\varphi(y^{1},y^{2})
=f\left(e^{\frac{h}{2}(y^{1}-y^{2})}-e^{\frac{h}{2}(y^{2}-y^{1})}, y^{1}+y^{2}\right)\in\O_{\CC^{2},p}^{W}
\end{equation*}
for $f(s^{1},s^{2})\in\O_{\CC^{2},\varphi(p)}$. 

We construct the inverse map $\psi_{p}:\O_{\CC^{2},p}^{W}\to\O_{\CC^{2},\varphi(p)}$ of $\varphi_{p}$. 
For every $n_{1},n_{2}\in\ZZ_{\ge 0}$, set the differential $D_{n_{1},n_{2}}:\O_{\CC^{2}}\to\O_{\CC^{2}}$ as
\begin{equation*}
D_{n_{1},n_{2}}:=\sum_{i=0}^{n_{1}+n_{2}}\left( \dfrac{(n_{1}+n_{2})!}{i!\,(n_{1}+n_{2}-i)!}\dfrac{1}{2^{n_{2}}h^{n_{1}}\left(e^{\frac{h}{2}(y^{1}-y^{2})}
+e^{\frac{h}{2}(y^{2}-y^{1})}\right)^{n_{1}}}\right)\left( \dfrac{\p}{\p y^{1}} \right)^{i}\left( \dfrac{\p}{\p y^{2}} \right)^{n_{1}+n_{2}-i} .
\end{equation*}
Then $\psi_{p}$ is given by 
\begin{equation*}
\psi_{p}(g)(s^{1},s^{2}):=\sum_{n_{1},n_{2}\in\ZZ_{\ge 0}}\dfrac{1}{n_{1}!\, n_{2}!}\left. D_{n_{1},n_{2}}(g)
\right|_{(y^{1},y^{2})=p}(s^{1}-s^{1}(p))^{n_{1}}(s^{2}-s^{2}(p))^{n_{2}}
\end{equation*}
for $g(y^{1},y^{2})\in\O_{\CC^{2},p}^{W}$. 
We have
\begin{align*}
& \left( \dfrac{\p}{\p s^{1}} \right)^{n_{1}}\left( \dfrac{\p}{\p s^{2}}\right)^{n_{2}} \\
=& \left( \dfrac{\p y^{1}}{\p s^{1}}\dfrac{\p}{\p y^{1}} +\dfrac{\p y^{2}}{\p s^{1}}\dfrac{\p}{\p y^{2}} \right)^{n_{1}}
\left( \dfrac{\p y^{1}}{\p s^{2}}\dfrac{\p}{\p y^{1}} +\dfrac{\p y^{2}}{\p s^{2}}\dfrac{\p}{\p y^{2}}\right)^{n_{2}} \\
=&\sum_{i=0}^{n_{1}+n_{2}}\left( \dfrac{(n_{1}+n_{2})!}{i!\,(n_{1}+n_{2}-i)!}\dfrac{1}{2^{n_{2}}h^{n_{1}}
\left(e^{\frac{h}{2}(y^{1}-y^{2})}+e^{\frac{h}{2}(y^{2}-y^{1})}\right)^{n_{1}}}\right)
\left( \dfrac{\p}{\p y^{1}} \right)^{i}\left( \dfrac{\p}{\p y^{2}} \right)^{n_{1}+n_{2}-i} \\
=&D_{n_{1},n_{2}}.
\end{align*}
The value of $\left. D_{n_{1},n_{2}}(g)\right|_{(y^{1},y^{2})=p}$ is determined uniquely by $(s^{1},s^{2})$ 
since $g$ is a $W$-invariant function.
Hence, we have $\psi_{p}=\varphi_{p}^{-1}$ by the Taylor expansion.
\qed
\end{pf}

By Definition \ref{analytic spaces} and Proposition \ref{prop:invariant theory of CC}, 
we have the following commutative diagram of complex manifolds:
\[
\xymatrix{
\widetilde{X}~ \ar@{>>}[d] \ar@{^{(}->}[r]& ~ \CC^{2} \ar@{>>}[d] \\
\widetilde{X}\ds W ~ \ar@{^{(}->}[r] & ~ \CC^{2}\ds W.
} 
\]

\begin{prop}
\label{prop:space of X ds W}
Set 
\begin{equation*}
\EE := \displaystyle \left\{z\in \CC ~\middle|~\frac{({\rm Re}\,z)^{2}}{\left(\exp\left(\frac{\pi^{2}}{\log \rho}\right)
-\exp\left(-\frac{\pi^{2}}{\log \rho}\right)\right)^{2}}+\frac{({\rm Im}\,z)^{2}}{\left(\exp\left(\frac{\pi^{2}}{\log \rho}\right)
+\exp\left(-\frac{\pi^{2}}{\log \rho}\right)\right)^{2}} < 1\right \} \subset \CC .
\end{equation*}
Then the morphism given by \eqref{eqn : isomorphism from C2 ds W to C2} induces an isomorphism between complex manifolds:
\begin{equation*}
\widetilde{X}\ds W \cong \EE\times \CC,\quad [(y^{1},y^{2})]\mapsto (s^{1},s^{2}).
\end{equation*}
In particular, the complex manifold $\widetilde{X}\ds W$ is isomorphic to the space of stability conditions ${\rm Stab}(\D^{b}(K_{\ell}))$ 
of $\D^{b}(K_{\ell})$. 
\end{prop}

\begin{pf}
It is easily shown that the complex manifold $\EE\times \CC$ is the image of $\widetilde{X}\ds W$ by $\underline{\varphi}$. 
Latter statement follows from ${\rm Stab}(\D^{b}(K_{\ell}))\cong \HH\times\CC$ in {\cite[Theorem 1.5]{dk}} and 
an isomorphism $\HH\cong \EE$ as complex manifolds. 
\qed
\end{pf}

\begin{rem}
The functions $e^{s^{1}}$ and $s^{2}$ are single-valued functions on $X$:
\begin{align*}
s^{1} &= \left( \dfrac{x^{1}}{x^{2}} \right)^{\frac{h}{2}}-\left( \dfrac{x^{2}}{x^{1}} \right)^{\frac{h}{2}},\\
e^{s^{2}} &= x^{1}x^{2}.
\end{align*}
Note that $s^{1}$ is an infinitely multi-valued function on $X$.
\end{rem}

\section{A Frobenius structure for $\ell$--Kronecker quiver}
Set $h:=\dfrac{2\pi\sqrt{-1}}{\log\rho}$. 
The eigenvalues $\rho$ and $\rho^{-1}$ of the Coxeter transformation $c$ can be expressed as follows:
\begin{equation*}
\rho=\exp\left(\dfrac{2\pi\sqrt{-1}}{h}\right),\quad\rho^{-1}=\exp\left(-\dfrac{2\pi\sqrt{-1}}{h}\right).
\end{equation*}
The value $h$ can be regarded as a generalization of the Coxeter number. 

We shall construct a Frobenius structure of rank $2$ and dimension $\displaystyle 1-\frac{2}{h}$ 
on the complex manifold $\widetilde{X} \ds W$ whose intersection form coincides with the generalized Cartan matrix $A$ for $K_{\ell}$. 
\begin{thm}
\label{thm : Frobenius manifold for l--Kronecker quiver}
There exists a unique Frobenius structure of rank $2$ and dimension $1-\dfrac{2}{h}$ on $\widetilde{X}\ds W$ 
such that $\displaystyle e=\frac{\p}{\p t^{1}}, \ E=t^{1}\dfrac{\p}{\p t^{1}}+\dfrac{2}{h}t^{2}\dfrac{\p}{\p t^{2}}$
for the flat coordinates $(t^{1}, t^{2})$ in Subsection \ref{flat coords 4.1} and 
the intersection form coincides with the generalized Cartan matrix $A$ of $\ell$-Kronecker quiver $K_{\ell}$. 
\end{thm}
\begin{rem}
Since the Frobenius structure is of rank $2$ and dimension $1-\dfrac{2}{h}$, the Euler vector field $E$ is given by
\begin{equation*}
E=t^{1}\dfrac{\p}{\p t^{1}}+\dfrac{2}{h}t^{2}\dfrac{\p}{\p t^{2}}
\end{equation*}
with respect to the flat coordinates $(t^{1},t^{2})$. In particular,  
it is easily shown that the Frobenius potential $\F$ must be
\begin{equation*}
\F(t^1,t^2)= \frac{1}{2}\eta_{12}(t^1)^2t^2+c(t^2)^{h+1}, \quad c\in \CC\setminus \{0\},
\end{equation*}
and hence the Frobenius structure on $\widetilde{X}\ds W$ should be unique if it exists (see also Example \ref{important example}).
\end{rem}

\subsection{The unit vector field and the Euler vector field}\label{flat coords 4.1}
Let $(s^{1},s^{2})$ be the coordinates of $\widetilde{X}\ds W$ in Proposition \ref{prop:space of X ds W}. 
Set
\begin{equation*}
t^{1}:=e^{\frac{h}{2}s^{2}}\cdot s^{1},\quad t^{2}:=e^{s^{2}}.
\end{equation*}
The map 
\begin{equation*}
\CC \times \EE \longrightarrow \CC\times\CC^{*};\quad (s^{1},s^{2})\mapsto(t^{1},t^{2})
\end{equation*}
is a local homeomorphism. 
The coordinates $(t^1,t^2)$ are locally equal to 
\begin{align*}
t^1=&~(x^{1})^{h}-(x^{2})^{h}, \\
t^2=&~x^{1}x^{2}
\end{align*}
where $(x^{1},x^{2})$ are the coordinates of $X$ in Definition \ref{coordinate change}. 
Define degrees of $t^{1}$ and $t^{2}$ as $\deg t^{1}:=h$ and $\deg t^{2}:=2$ due to Proposition \ref{degree two coordinate} and 
local expressions of $t_{1}$ and $t_{2}$ above.
Set 
\begin{equation*}
e:=\dfrac{\p}{\p t^{1}},\quad E:=\dfrac{\deg t^{1}}{h}t^{1}\dfrac{\p}{\p t^{1}}+\dfrac{\deg t^{2}}{h}t^{2}\dfrac{\p}{\p t^{2}}=t^{1}\dfrac{\p}{\p t^{1}}+\dfrac{2}{h}t^{2}\dfrac{\p}{\p t^{2}} .
\end{equation*}
It will be shown later that $(t^{1},t^{2})$ are flat coordinates and $e$ and $E$ are the unit vector field and the Euler vector field. 
\begin{rem}
Note that the Euler vector field $E$ is equal to $\dfrac{2}{h}\dfrac{\p}{\p s^{2}}$. This observation can be regarded as an analogy of
\cite[Theorem 2.1 (ii)]{dz}.
\end{rem}

\subsection{The non-degenerate $\O_{\widetilde{X}\ds W}$-symmetric bilinear form and flat coordinates}
Let $g:\Omega_{X}^{1}\times\Omega_{X}^{1}\to\O_{X}$ be the bilinear form on the cotangent sheaf of $X$ induced by 
the generalized Cartan matrix of $K_{\ell}$ via the natural isomorphism $T^{*}_{p}\h\cong \h^{*}$ and the restriction to $X$:
\begin{equation*}
g(dx^i,dx^j):=I(\beta_i^*,\beta_j^*)=
\begin{pmatrix}
0 & -1 \\
-1 & 0
\end{pmatrix}.
\end{equation*}
The bilinear form $g$ induces the one on the cotangent sheaf of $\widetilde{X}$ 
(denote the induced bilinear form on $\Omega_{\widetilde{X}}^{1}$ by the same symbol $g$ for simplicity). 
With respect to the coordinates $(y^{1},y^{2})$ of $\widetilde{X}$ in \eqref{defn:universal covering  of X}, we have 
\begin{align*}
g(dy^i,dy^j)&=\sum_{a,b=1}^2\dfrac{\p y^i}{\p x^a}\dfrac{\p y^j}{\p x^b} g(dx^a,dx^b) \\
&=
\begin{pmatrix}
0  & -\dfrac{1}{e^{y^{1}+y^{2}}} \\
-\dfrac{1}{e^{y^{1}+y^{2}}} & 0
\end{pmatrix}.
\end{align*}
The bilinear form $g$ on $\Omega_{\widetilde{X}}^{1}$ induces the bilinear form on $\Omega_{\widetilde{X}\ds W}^{1}$
(denote this induced bilinear form by the same symbol $g$ again). 
By the definition of the coordinates $(s^{1},s^{2})$ in Proposition \ref{prop:invariant theory of CC}, the induced bilinear form
$g:\Omega_{\widetilde{X}\ds W}^{1}\times\Omega_{\widetilde{X}\ds W}^{1}\to\O_{\widetilde{X}\ds W}$ is given by 
\begin{align}
\nonumber
g(ds^{i},ds^{j})&=\sum_{a,b=1}^{2}\dfrac{\p s^{i}}{\p y^{a}}\dfrac{\p s^{j}}{\p y^{b}} g(dy^{a},dy^{b}) \\
&\label{eqn : metric on widetilde X ds W}=
\begin{pmatrix}
-\dfrac{2}{e^{s^{2}}}&0\\
0&\dfrac{h^{2}}{2e^{s^{2}}}(4+(s^{1})^{2})
\end{pmatrix}.
\end{align}
\begin{prop}\label{desired bilinear form}
Define $\eta:\Omega^{1}_{\widetilde{X}\ds W}\times\Omega^{1}_{\widetilde{X}\ds W}\to\O_{\widetilde{X}\ds W}$ as 
\begin{equation*}
\eta:={\rm Lie}_{e}\,g \ .
\end{equation*}
Then the bilinear form $\eta$ defines a non-degenerate and flat $\O_{\widetilde{X}\ds W}$--symmetric bilinear form on 
$\T_{\widetilde{X}\ds W}$ (denote this bilinear form on $\T_{\widetilde{X}\ds W}$ by the the same symbol $\eta$). 
\end{prop}
\begin{pf}
With respect to the local coordinates $(t^{1},t^{2})$, we have
\begin{align*}
\eta\left( dt^{i},dt^{j} \right)&=\dfrac{\p}{\p t^{1}}\left( g\left( dt^{i},dt^{j}\right) \right)
-g\left( {\rm Lie}_{\frac{\p}{\p t^{1}}}dt^{i},dt^{j}\right)-g\left( dt^{i},{\rm Lie}_{\frac{\p}{\p t^{1}}}dt^{j}\right) \\
&=\dfrac{\p}{\p t^{1}}\left( g\left( dt^{i},dt^{j}\right) \right) .
\end{align*}
By the equation \eqref{eqn : metric on widetilde X ds W}, we have 
\begin{equation*}
 (g(dt^{i} , dt^{j} )) = \left(
\begin{array}{cc}
2h^{2} ( t^{2})^{h-1} & -h t^{1} \\
-h t^{1} & -2 t^{2}
\end{array}
\right) 
\end{equation*}
and hence
\begin{equation*}
\displaystyle
(\eta(dt^{i} , dt^{j} )) = \left(
\begin{array}{cc}
0&-h\\
-h&0
\end{array}
\right), 
\ \ \text{i.e.,} \ \
\left(\eta\left(\frac{\p}{\p t^{i}}, \frac{\p}{\p t^{j}}\right)\right)=\left(
\begin{array}{cc}
0&-\frac{1}{h}\\
-\frac{1}{h}&0
\end{array}
\right).
\end{equation*}
Therefore, the bilinear form $\eta$ is a non-degenerate and flat bilinear form whose flat coordinates are $(t^{1},t^{2})$.
\qed
\end{pf}

\subsection{The product structure}\label{reconst of product from g}
Recall Lemma~\ref{lem : dubrovin}.
Set $\Gamma^{ij}_{k}:=g\left( dt^{i},\nabla_{\frac{\p}{\p t^{k}}} dt^{j} \right)$ 
where $\nabla$ is the Levi--Civita connection of $g$ on $\Omega^{1}_{\widetilde{X} \ds W}$.
Define $C_{ij}^{k}$ as follows:
\begin{equation*}
C_{ij}^{k}:=\dfrac{h}{\deg t^{k}-1}\sum_{a=1}^{2}\eta_{ia}\cdot \Gamma^{ak}_{j},
\quad i,j,k\in \{1,2\}
\end{equation*}
where $\displaystyle \eta_{ij}:=\eta\left( \frac{\p}{\p t^{i}},\frac{\p}{\p t^{j}} \right)$ in Proposition \ref{desired bilinear form}. 
\begin{lem}\label{prop:product and potential}
Set
\begin{equation*}
\displaystyle \F:=-\frac{1}{2h}\left( t^{1} \right)^{2}t^{2}+\frac{1}{h^{2}-1}\left( t^{2} \right)^{h+1}.
\end{equation*}
Then we have 
\begin{equation*}
C_{ij}^{k}=\sum_{a=1}^{2}\eta^{ka}\dfrac{\p^{3}\F}{\p t^{a}\p t^{i}\p t^{j}}\ .
\end{equation*}
\end{lem}
\begin{pf}
The statement follows from straightforward calculations. 
\qed
\end{pf}

The following proposition follows from Lemma \ref{prop:product and potential}: 
\begin{prop}
Define the $\O_{\widetilde{X}\ds W}$-linear map
$\circ: \T_{\widetilde{X}\ds W}\times \T_{\widetilde{X}\ds W}\rightarrow \T_{\widetilde{X}\ds W}$ by 
\begin{equation*}
\displaystyle \frac{\p}{\p t^{i}}\circ\frac{\p}{\p t^{j}}:=\sum_{k=1}^{2}C_{ij}^{k}\frac{\p}{\p t^{k}},\quad i,j=1,2.
\end{equation*}
Then the product $\circ$ is commutative and associative. 
Moreover, $e=\dfrac{\p}{\p t^{1}}$ is the unit vector field for $\circ$. 
\end{prop}

\subsection{The discriminant locus and the monodromy group}
Denote by $\D_{\widetilde{X}\ds W}$ the zeros of the discriminant $\Delta_{\widetilde{X}\ds W}$ 
of the Frobenius manifold $({\widetilde{X}\ds W},\eta,\circ,e,E)$. 
The discriminant locus $\D_{\widetilde{X}\ds W}$ is given by 
\begin{equation*}
\D_{\widetilde{X}\ds W}=\left\{ (s^{1},s^{2})\in\widetilde{X}\ds W ~\middle|~ s^{1}=\pm2\sqrt{-1} \right\}
\end{equation*}
with respect to the global coordinates $(s^{1},s^{2})$ by the equation \eqref{eqn : metric on widetilde X ds W}. 
With respect to the flat coordinates $(t^{1},t^{2})$, we have
\begin{equation*}
\D_{\widetilde{X}\ds W}=\left\{ (t^{1},t^{2})\in\widetilde{X}\ds W ~\middle|~ (t^{1})^{2}+4(t^{2})^{h}=0 \right\} .
\end{equation*}

\begin{prop}
Set 
\[
\widetilde{X}_{\rm reg}:=\pi_{X}^{-1}(X_{\rm reg})\subset \widetilde{X}.
\]
Then we have 
\begin{equation}
\label{eqn : discriminant}
\widetilde{X}_{\rm reg}/W \cong\left(\widetilde{X}\ds W\right)\backslash\D_{\widetilde{X}\ds W}.
\end{equation}
\end{prop}
\begin{pf}
Recall Proposition \ref{proper disconti and free}.
Since $\bigcup_{\alpha\in\Delta^{\rm re}_+}H_\alpha=\{(x^1,x^2)~|~x^2=\nu^{2k+1}x^1~(k\in\ZZ)\}\subset X$, 
we have
\begin{equation*}
\pi_X^{-1}\left(\bigcup_{\alpha\in\Delta^{\rm re}_+}H_\alpha\right)=
\left\{(y^{1},y^{2})\in\widetilde{X}~|~y^{1}=y^{2}+(2k+1)\log\nu,~k\in\ZZ\right\}\subset \widetilde{X}
\end{equation*}
and hence 
\begin{equation*}
\widetilde{X}\backslash\pi_{X}^{-1}\left(\bigcup_{\alpha\in\Delta^{\rm re}_+}H_\alpha\right)=\widetilde{X}_{\rm reg}.
\end{equation*}
By the definition of $(s^1,s^2)$, the image of 
$(y^{1},y^{2})\in\pi_{X}^{-1}\left(\bigcup_{\alpha\in\Delta^{\rm re}_+}H_\alpha\right)$ by $\underline{\varphi}$ are 
\begin{align*}
s^{1} & = \exp \left( - \dfrac{h}{2} (2k+1) \log \nu \right) - \exp \left( \dfrac{h}{2} (2k+1) \log \nu \right) \\
& = \exp \left( - \left( k + \dfrac{1}{2} \right) \pi \sqrt{-1} \right) - \exp \left( \left( k + \dfrac{1}{2} \right) \pi \sqrt{-1} \right) \\
& = \left\{
\begin{array}{cl}
-2\sqrt{-1} & ( \text{if $k$ is even}  ) \\
2\sqrt{-1} & ( \text{if $k$ is odd} ),
\end{array}
\right.\\
s^{2} & = 2 y^{1} + (2k+1) \log\nu.
\end{align*}
Therefore the isomorphism \eqref{eqn : discriminant} holds. 
\qed
\end{pf}

\begin{prop}\label{monodromy prop}
The monodromy group $W_{\widetilde{X}\ds W}$ of $(\widetilde{X}\ds W,\eta,\circ,e,E)$ coincides with the Weyl group $W$.
\end{prop}

\begin{pf} 
By the construction, the coordinates $x^{1}$ and $x^{2}$ are periods of the Frobenius manifold $(\widetilde{X}\ds W,\eta,\circ,e,E)$: 
\begin{equation*}
\begin{array}{cl}
x^{1} &=e^{\frac{s^{2}}{2}}\left( \dfrac{s^{1}+\sqrt{\left(s^{1}\right)^{2}+4}}{2} \right)^{\frac{1}{h}}, \\
x^{2} &=e^{\frac{s^{2}}{2}}\left( \dfrac{-s^{1}+\sqrt{\left(s^{1}\right)^{2}+4}}{2} \right)^{\frac{1}{h}}.
\end{array}
\end{equation*}
Fix a point $*\in\widetilde{X}_{\rm reg}$. 
For $i=1,2$, let $\widetilde{\gamma_{i}}$ be a path on $\widetilde{X}$ from $*$ to $r_{i}(*)$. 
Denote by $\gamma_{i}:=\underline{\varphi}\left(\widetilde{\gamma_{i}}\right)$. 
Then $\gamma_{1}$ (resp. $\gamma_{2}$) is a loop around $2\sqrt{-1}$ (resp. $-2\sqrt{-1}$) on $(\widetilde{X}\ds W)_{\rm reg}$ 
and these loops are generators of the fundamental group 
$\pi_{1}\left( ( \widetilde{X}\ds W)_{\rm reg} , \underline{\varphi}(*) \right)$.
Note that the branch cut of $\sqrt{\left(s^{1}\right)^{2}+4}$ is the line between $2\sqrt{-1}$ and $-2\sqrt{-1}$. 
The monodromy representation $M_{\gamma_{i}}$ of $\gamma_{i}$ is given by 
\[
M_{\gamma_{1}}=
\begin{pmatrix}
0 & \nu \\
\nu^{-1} & 0
\end{pmatrix},
\quad 
M_{\gamma_{2}}=
\begin{pmatrix}
0 & \nu^{-1} \\
\nu & 0
\end{pmatrix}
\]
for periods $(x^{1},x^{2})$. 
Therefore, we have $W_{\widetilde{X}\ds W}\cong W$.
\qed 
\end{pf}

Finally, we show that the Frobenius manifold $(\widetilde{X}\ds W,\eta,\circ,e,E)$ is semi-simple
in the sense of {\cite[Lecture 3]{dub}}:
\begin{prop}
The Frobenius manifold $(\widetilde{X}\ds W,\eta,\circ,e,E)$ is semi-simple. 
\end{prop}
\begin{pf}
Set 
\begin{equation*}
u^{1}:=e^{h s^{2}}\left(s^{1}+2\sqrt{-1}\right),\quad u^{2}:=e^{h s^{2}}\left(s^{1}-2\sqrt{-1}\right).
\end{equation*}
We can easily check that $(u^{1},u^{2})$ form canonical coordinates. 
\qed
\end{pf}
\subsection{Perspectives}

In the case of the root system of type $A_{2}$, there exsits an isomorphism between the Frobenius manifold $\h /W$ and $\Stab(\D^{b}(A_{2}))$ such that the central charge map is identified with the oscillatory integrals on $\h/W$ (\cite{bqs,hkk}). 
Based on the results and Proposition \ref{prop:space of X ds W}, we expect the following 
\begin{conj}
\label{conj : Frobenius manifold and satbility conditions for generalized Kronecker quiver}
There should exist an isomorphism $\varphi : \widetilde{X}\ds W \to \Stab( \D^b ( K_{\ell} ))$ such that $\varphi$ is compatible with a deformed flat coordinates (see \cite{dub3}) and the central charge map. 
That is, the following diagram commutes; 
\begin{equation*}
\xymatrix{
\widetilde{X}\ds W \ar[rr]^\varphi \ar[rd]_{( \widetilde{t}^{1}|_{u=1} , \widetilde{t}^{2}|_{u=1})} & & \Stab( \D^b (K_\ell )) \ar[ld]^\Z \\
& \CC^2 &
} .
\end{equation*}
In particular, we have 
\begin{equation*}
\widehat{\nabla}|_{u=1} = \varphi^{*} d , 
\end{equation*}
where $d$ is the trivial connection on 
\begin{equation*}
T \Stab(\D^{b}(K_{\ell})) \cong \Stab(\D^{b}(K_{\ell})) \times \Hom_{\ZZ}(K_{0}(\D^{b}(K_{\ell})),\CC)
\end{equation*}
and $\widehat{\nabla}|_{u=1}$ is the restriction of the first structure connection on  $T(\widetilde{X}\ds W)$, which is defined by 
\begin{equation*}
\widehat{\nabla}_{\delta}\delta' := \ns_{\delta'}\delta+\dfrac{1}{u}\delta\circ\delta' .
\end{equation*}
\end{conj}

Through the present work, we reach to partial results that exponents, duality among them and 
the non-degenerate symmetric bilinear form can be obtained via eigenvalues and eigenvectors of 
the Coxeter transformation for a generalized root system whose Cartan matrix is non-degenerate and symmetric
(e.g., see Lemma~\ref{coordinate change}).  

The present paper will be the first of a series of our attempts to construct Frobenius structures from 
the Weyl group invariant theories associated to such generalized root systems, moreover on the spaces of stability conditions. 

\end{document}